\documentclass[12pt]{amsart}
\usepackage{graphicx}
\usepackage{amsfonts}
\usepackage{amsmath}
\usepackage{amssymb}
\newtheorem{thm}{Theorem}[section]
\newtheorem{lem}[thm]{Lemma}
\newtheorem{cor}[thm]{Corollary}
\newtheorem{prp}[thm]{Proposition}

\theoremstyle{definition}

\theoremstyle{remark}
\newtheorem{remark}[thm]{Remark}
\numberwithin{equation}{section}
\DeclareMathOperator{\com}{Comp}
\DeclareMathOperator{\co}{\mathcal{C}}
\DeclareMathOperator{\col}{Col}
\DeclareMathOperator{\Des}{Des}
\DeclareMathOperator{\C}{C}

\title{The $m$-colored composition poset}

\author[B. Drake]{Brian Drake}
\address{Department of Mathematics, Brandeis University, Waltham, MA, USA, 02454}
\email{bdrake@brandeis.edu}
\urladdr{http://people.brandeis.edu/\~{}bdrake}

\author[T. K. Petersen]{T. Kyle Petersen}
\address{Department of Mathematics, Brandeis University, Waltham, MA, USA, 02454}
\email{tkpeters@brandeis.edu}
\urladdr{http://people.brandeis.edu/\~{}tkpeters}


\begin{document}
\begin{abstract}
We generalize Bj\"{o}rner and Stanley's poset of compositions to $m$-colored compositions. Their work draws many analogies between their (1-colored) composition poset and Young's lattice of partitions, including links to (quasi-)symmetric functions and representation theory. Here we show that many of these analogies hold for any number of colors. While many of the proofs for Bj\"{o}rner and Stanley's poset were simplified by showing isomorphism with the subword order, we remark that with 2 or more colors, our posets are not isomorphic to a subword order.
\end{abstract}

\maketitle

\section{Introduction}

This paper explores a generalization of the poset of compositions introduced in recent work of Bj\"{o}rner and Stanley \cite{BjornerStanley}, which draws several analogies between their poset and Young's lattice of partitions. We recall some key facts about Young's lattice.

A partition $\lambda = (\lambda_1, \lambda_2, \ldots )$ of $n$, denoted $\lambda \vdash n$, is a sequence of nonnegative integers $\lambda_1 \geq \lambda_2 \geq \cdots \geq 0$ such that $\sum \lambda_i = n$. The set of all partitions of all integers $n \geq 0$ forms a lattice under the partial order given by inclusion of Young diagrams: $\lambda \leq \mu$ if $\lambda_i \leq \mu_i$ for all $i$. This lattice is called \emph{Young's lattice}, $Y$, which has several remarkable properties, including the following list given in \cite{BjornerStanley}.
\begin{enumerate}
\item[Y1.] $Y$ is a graded poset, where a partition $\lambda \vdash n$ has rank $n$.

\item[Y2.] The number of saturated chains from the minimal partition $\emptyset$ to $\lambda$ is the number $f^{\lambda}$ of Young tableaux of shape $\lambda$.

\item[Y3.] The number of saturated chains from $\emptyset$ to rank $n$ is the number of involutions in the symmetric group $\mathfrak{S}_n$.

\item[Y4.] Let $s_{\lambda}$ denote a Schur function. Pieri's rule \cite{Stanley2} gives \[ s_1s_{\lambda} = \sum_{\lambda \prec \mu} s_{\mu}, \] where $\lambda \prec \mu$ means that $\mu$ covers $\lambda$ in $Y$.

\item[Y5.] Since $Y$ is in fact a distributive lattice, every interval $[\lambda,\mu]$ is EL-shellable and hence Cohen-Macaulay.

\item[Y6.] $Y$ is the Bratteli diagram for the tower of algebras $K \mathfrak{S}_0 \subset K\mathfrak{S}_1 \subset \cdots$, where $K\mathfrak{S}_n$ denotes the group algebra of $\mathfrak{S}_n$ over $K$, a field of characteristic zero.
\end{enumerate}

For each of the properties listed above, Bj\"{o}rner and Stanley
give an analogous property for their poset of compositions. Here we
will consider more generally the poset of $m$-colored compositions,
and show that this poset has properties analogous to Y1--5 above.
Finding an analog of property Y6 is an open problem.

Recall that a composition $\alpha = (\alpha_1, \alpha_2, \ldots,
\alpha_k)$ is an ordered tuple of positive integers, called the \emph{parts} of $\alpha$. We write $l(\alpha) = k$ for the number of parts of $\alpha$. If the sum of the parts of $\alpha$ is $n$, i.e., $|\alpha| := \alpha_1 + \alpha_2 + \cdots + \alpha_k = n$, then we say $\alpha$ is a composition of $n$, written $\alpha \models n$. Let $\com(n)$ denote
all the compositions of $n$, and define the set of all compositions
\[\co := \bigcup_{n \geq 0} \com(n),\] where $\emptyset$ is the unique
composition of $0$.

Bj\"{o}rner and Stanley give $\co$ a partial order defined by the
following covering relations. We say $\beta$ \emph{covers} $\alpha$, written $\alpha \prec \beta$, if $\alpha < \beta$ and there is no $\beta'$ such that $\alpha < \beta' < \beta$. We can obtain all compositions $\beta$ that cover $\alpha$ by adding 1 to a part of $\alpha$ or adding 1 to a part and splitting that part
into two parts. In other words, for some $j$ we can write $\beta$ as:
\begin{enumerate}
\item $( \alpha_1, \ldots, \alpha_{j-1}, \alpha_{j}+1, \alpha_{j+1}, \ldots, \alpha_k)$, or

\item $( \alpha_1, \ldots, \alpha_{j-1}, h+1, \alpha_{j}-h, \alpha_{j+1}, \ldots, \alpha_k)$ for some $0 \leq h \leq \alpha_j-1$.
\end{enumerate}

We can generalize this poset to a family of posets indexed by
the number of ``colors" $m \geq 1$. An $m$-colored composition is an
ordered tuple of colored positive integers, say $\alpha = (
\varepsilon_1 \alpha_1, \varepsilon_2 \alpha_2, \ldots,
\varepsilon_k \alpha_k)$, where the $\alpha_s$ are positive integers and if $\omega$ is a primitive $m$-th
root of unity, $\varepsilon_s = \omega^{i_s}$, $0 \leq i_s \leq
m-1$. We say the part $\varepsilon_s \alpha_s$ has color $\varepsilon_s$, and we write
$\alpha \models_m n$ if $|\alpha|:= \alpha_1 + \alpha_2 + \cdots + \alpha_k =
n$. For example, if $m =3$, then $\alpha = ( \omega 2, 1, \omega^2 1, 3 )$ is a 3-colored composition of $2+1+1+3 = 7$.

Note that there are $m^k$ ways to color any ordinary composition
of $n$ with $k$ parts, leading us to conclude that there are \[
\sum_{k=1}^n \binom{n-1}{k-1} m^k = m(m+1)^{n-1} \] $m$-colored
compositions of $n$ (so if $m=1$, we have $2^{n-1}$ ordinary compositions). Let $\com^{(m)}(n)$ denote the set of all
$m$-colored compositions of $n$, and define \[\co^{(m)} := \bigcup_{n
\geq 0} \com^{(m)}(n),\] where $\emptyset$ is again the unique
composition of $0$.

\begin{figure}
\centering
\includegraphics[angle = 90, scale = .9]{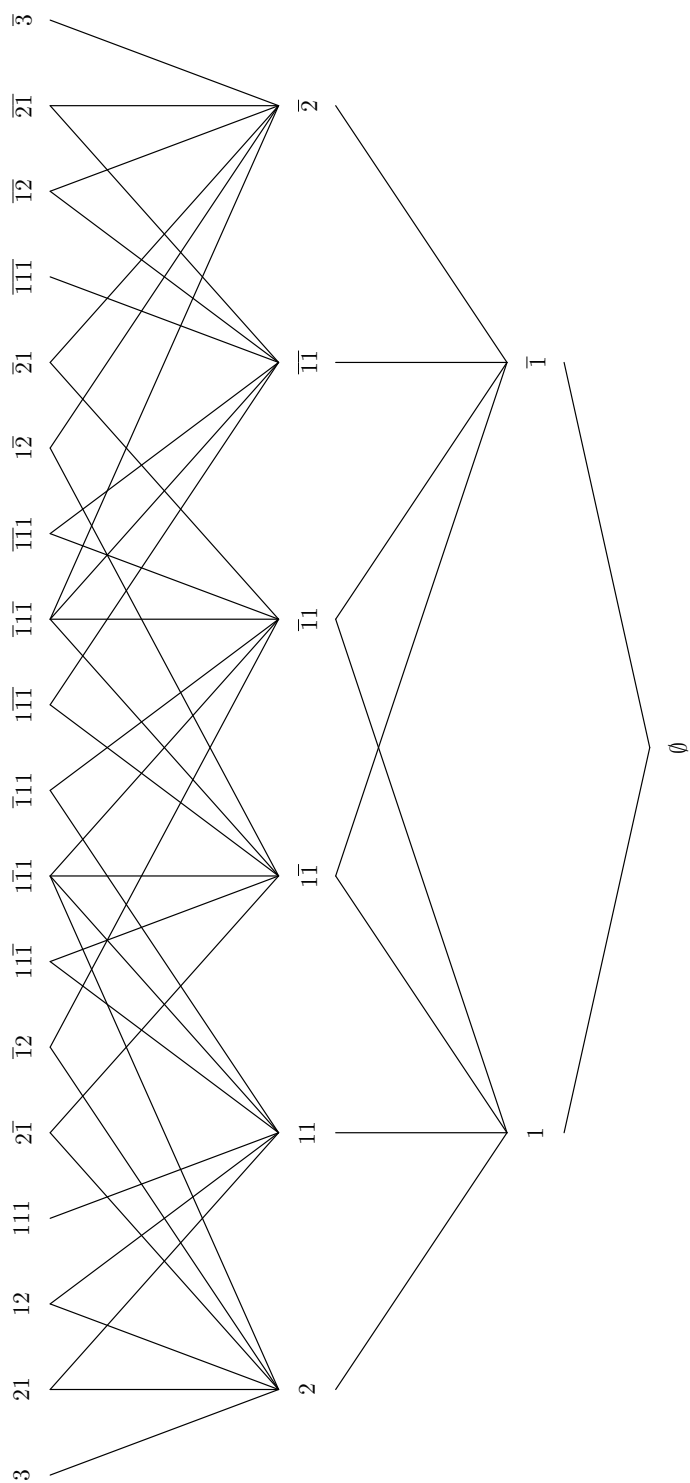}
\caption{The first four levels of the 2-colored composition poset.\label{fig:2color}}
\end{figure}

We can define a partial order on $\co^{(m)}$ with many of the same
properties of $\co$. The covering relations are as follows. We have
$\beta$ covers $\alpha$ if, for some $j$, we can write $\beta$ as:
\begin{enumerate}
\item $( \varepsilon_1 \alpha_1, \ldots, \varepsilon_{j-1} \alpha_{j-1}, \varepsilon_j(\alpha_{j}+1), \varepsilon_{j+1} \alpha_{j+1}, \ldots, \varepsilon_k \alpha_k )$,

\item $( \varepsilon_1 \alpha_1, \ldots, \varepsilon_{j-1} \alpha_{j-1}, \varepsilon_j (h+1), \varepsilon_j(\alpha_{j}-h), \varepsilon_{j+1} \alpha_{j+1}, \ldots, \varepsilon_k \alpha_k )$ for some $0 \leq h \leq \alpha_j-1$, or

\item $( \varepsilon_1 \alpha_1, \ldots, \varepsilon_{j-1} \alpha_{j-1}, \varepsilon_j h, \varepsilon' 1,  \varepsilon_j (\alpha_{j}-h), \varepsilon_{j+1} \alpha_{j+1}, \ldots, \varepsilon_k \alpha_k )$ where $\varepsilon' \neq \varepsilon_j$ and $0 \leq h \leq \alpha_j-1$, with the understanding that we will ignore parts of size 0.
\end{enumerate}
Relations (1) and (2) are just like those of $\co$: while preserving the color, we add 1 to a part, or we add 1 to a part and split that part into two parts. Relation (3) handles the case where the color of the ``1" we add differs from where we try to add it. Notice that it is immediate from these cover relations that $\co^{(m)}$ is a graded poset with level $n$ consisting of all $m$-compositions of $n$. This property is analogous to property Y1 of Young's lattice. See Figure \ref{fig:2color} for the first four levels of the 2-colored composition poset.

In this paper we will show that for any positive fixed $m$, the
poset $\co^{(m)}$ possesses properties analogous to $\co =
\co^{(1)}$ (and indeed to Young's lattice). In fact many of the
arguments used in \cite{BjornerStanley} generalize in a
straightforward way. One important argument that doesn't generalize is the isomorphism shown between $\co$ and the subword order; see Remark \ref{rem:1}. In section \ref{sec:descent} we discuss colored
permutations, their color-descent compositions, and chains in
$\co^{(m)}$.  In section \ref{sec:quasi} we present Poirier's colored
quasisymmetric functions \cite{Poirier} and show that $\co^{(m)}$ gives a
Pieri-type rule for multiplying a fundamental basis.  We define a
CL-labeling in section \ref{sec:mobius} and use this to calculate
the M\"{o}bius function of lower intervals.  Section
\ref{sec:clproof} contains the proof that this labeling is a
CL-labeling.

\section{Colored permutations and descent sets}\label{sec:descent}

Compositions can be used to encode descent classes of ordinary permutations in the following way. Recall that a \emph{descent} of a permutation $w \in \mathfrak{S}_n$ is a position $i$ such that $w_i > w_{i+1}$, and that an \emph{increasing run} (of length $r$) of a permutation $w$ is a maximal subword of consecutive letters $w_{i+1} w_{i+2} \cdots w_{i+r}$ such that $w_{i+1} < w_{i+2} < \cdots < w_{i+r}$. By maximality, we have that if $w_{i+1} w_{i+2} \cdots w_{i+r}$ is an increasing run, then $i$ is a descent of $w$ (if $i\neq 0$), and $i+r$ is a descent of $w$ (if $i+r \neq n$). For any permutation $w \in \mathfrak{S}_n$ define the \emph{descent composition}, $\C(w)$, to be the ordered tuple listing from left to right the lengths of the increasing runs of $w$. If $\C(w) = (\alpha_1, \alpha_2, \ldots, \alpha_k)$, we can recover the descent set of $w$:
\[ \Des(w) := \{ i\, : w_i > w_{i+1} \} = \{ \alpha_1, \alpha_1 + \alpha_2, \ldots, \alpha_1 + \alpha_2 + \cdots + \alpha_{k-1} \}.\] For example, the permutation $w = 345261$ has $\C(w) = (3,2,1)$ and $\Des(w) = \{ 3, 5\}$. We now define colored permutations and colored descent compositions.

Loosely speaking, $m$-colored permutations are permutations where
each of the elements permuted are given one of $m$ ``colors." If $\omega$ is any primitive $m$-th root of unity, \[
\omega 3 \,\, \omega 2 \,\, 1 \,\, \omega^3 4 \] is an example of a
colored permutation. We can think of building colored permutations
by taking an ordinary permutation and then arbitrarily assigning
colors to the letters, so we see that there are $m^n n!$ $m$-colored
permutations of $[n]:=\{1,2,\ldots,n\}$.

Strictly speaking, $m$-colored permutations are elements of the
wreath product $C_m \wr \mathfrak{S}_n$, where $C_m = \{ 1,\omega,\ldots, \omega^{m-1}\}$ is the cyclic
group of order $m$.  We write an element $u = u_1 u_2 \cdots u_n \in
C_m \wr \mathfrak{S}_n$ as a word in the alphabet \[C_m \times [n]:= \{ 1,2,\ldots,n, \omega 1, \omega 2,\ldots,\omega n, \ldots, \omega^{m-1} 1, \omega^{m-1} 2, \ldots, \omega^{m-1} n\}, \] such that $|u| = |u_1| |u_2|
\cdots |u_n|$ is an ordinary permutation in $\mathfrak{S}_n$. We say
$\varepsilon_i = u_i/|u_i|$ is the color of $u_i$.

For any $u \in C_m \wr \mathfrak{S}_n$, we can write $u = v_1 v_2 \cdots v_k$ so that each $v_i$ is a word in which all the letters have the same color, $\varepsilon'_i$, and no two consecutive colors are the same: $\varepsilon'_i \neq \varepsilon'_{i+1}$, $i = 1,2,\ldots, k-1$.
Then we define the \emph{color composition} of $u$, \[\col(u) := (
\varepsilon'_1 \alpha'_1, \varepsilon'_2 \alpha'_2, \ldots,
\varepsilon'_k \alpha'_k),\] where $\alpha'_s$ denotes the number of letters in $v_s$. Now suppose an $m$-colored permutation $u$ has color composition $\col(u) = $ $(\varepsilon_1' \alpha'_1, \varepsilon_2' \alpha'_2, \ldots,
\varepsilon_k' \alpha'_k)$.  Then the \emph{colored descent composition}
\[\C^{(m)}(u) := ( \varepsilon_1 \alpha_1, \varepsilon_2 \alpha_2, \ldots, \varepsilon_l \alpha_l), \]
is the refinement of $\col(u)$ where we replace part $\varepsilon_i'
\alpha_i'$ with $\varepsilon_i' \C(|v_i|)$, where $\C$ is the ordinary descent composition, and we view $|v_i|$ as an ordinary permutation of distinct letters.

More intuitively, the colored descent composition $\C^{(m)}(u)$ is the ordered tuple listing the lengths of increasing runs of $u$ with constant color, where we record not only the length of such a run, but also its color. An example should cement the notion. If we have two colors (indicated with a bar), let \[u = 1\bar{2}\bar{3}4\bar{8}\bar{5}76.\]
Then the color composition is $\col(u) = (1,\bar{2},1,\bar{2},2)$, and \[\C^{(m)}(u) = (1,\bar{2},1,\bar{1},\bar{1},1,1).\]

For any $\alpha \in \com^{(m)}(n)$, a \emph{saturated chain} from $\emptyset$ to $\alpha$ is a sequence of compositions \[ \emptyset = \alpha^0 \prec \alpha^1 \prec \cdots \prec \alpha^n = \alpha,\] where $\prec$ denotes a cover relation in $\co^{(m)}$, and therefore $\alpha^i \in \com^{(m)}(i)$. Now, given any $u \in C_m \wr \mathfrak{S}_n$, let $u[i]$ denote the restriction of $u$ to letters in $C_m \times [i]$. For example, if $u = \bar{2}17\bar{6}\bar{3}\bar{4}58$, then $u[5] = \bar{2}1\bar{3}\bar{4}5$. We then define the sequence \[\mathfrak{m}(u) := ( \C^{(m)}(u[1]), \ldots, \C^{(m)}(u[n])),\] so that $\C^{(m)}(u[i]) \in \com^{(m)}(i)$. Using the same example $u =\bar{2}17\bar{6}\bar{3}\bar{4}58$, we have \[ \mathfrak{m}(u) = ( 1, \bar{1}1, \bar{1}1\bar{1}, \bar{1}1\bar{2}, \bar{1}1\bar{2}1, \bar{1}1\bar{1}\bar{2}1, \bar{1}2\bar{1}\bar{2}1, \bar{1}2\bar{1}\bar{2}2).\] The following theorem is the natural generalization of Theorem 2.1 of \cite{BjornerStanley}.

\begin{thm}\label{thm:saturated}
The map $\mathfrak{m}$ is a bijection from $C_m \wr \mathfrak{S}_n$ to saturated chains from $\emptyset$ to $\alpha$, where $\alpha$ ranges over all colored compositions in $\com^{(m)}(n)$.
\end{thm}

This proof follows the same line of reasoning used by Bj\"{o}rner and Stanley in proving Theorem 2.1 of \cite{BjornerStanley}.

\begin{proof}
For any colored permutation $u \in C_m \wr \mathfrak{S}_n$ define, for all $0\leq i \leq n$ and all $0\leq j \leq m-1$, \[ u_{(i,j)} := u_1 \cdots u_i \, \omega^j(n+1) \, u_{i+1} \cdots u_n.\] In other words, the $u_{(i,j)}$ are all those permutations $w$ in $C_m \wr \mathfrak{S}_{n+1}$ such that $w[n] = u$. We will show that the compositions $\C^{(m)}(u_{(i,j)})$ are all distinct and moreover that they are precisely those compositions in $\com^{(m)}(n+1)$ that cover $\C^{(m)}(u)$.

Suppose $\C^{(m)}(u) = ( \varepsilon_1 \alpha_1, \ldots, \varepsilon_k \alpha_k)$, and let $b_s = \alpha_1 + \cdots + \alpha_s$, with the convention that $b_0 = 0$. For any fixed $j = 0,1,\ldots, m-1$, we have two cases, corresponding to cover relations of type (1) or type (3):
\[ \C^{(m)}(u_{(b_s,j)}) = \begin{cases} (\varepsilon_1 \alpha_1, \ldots, \varepsilon_s (\alpha_s + 1), \ldots, \varepsilon_k \alpha_k) &  \mbox{ if } \varepsilon_s = \omega^j, \\
(\varepsilon_1 \alpha_1, \ldots, \varepsilon_s \alpha_s, \omega^j 1, \ldots, \varepsilon_k \alpha_k) & \mbox{ otherwise.}
\end{cases}\]
All these compositions, over $s = 0,\ldots,k$, $j = 0,\ldots, m-1$, are distinct and cover $\C^{(m)}(u)$. To consider the other cases, suppose $i$ is not of the form $\alpha_1 + \cdots + \alpha_s$. Then it can be written as $i = \alpha_1 + \cdots + \alpha_s + h$, where $0 \leq s \leq k$ and $1 \leq h \leq \alpha_{s+1}-1$ (if $s = 0$, then $i = h$). Again we have two cases, corresponding to cover relations of type (2) or type (3):
\[ \C^{(m)}(u_{(i,j)}) = \begin{cases}(\varepsilon_1 \alpha_1, \ldots, \varepsilon_{s+1}(1+h), \varepsilon_{s+1} (\alpha_{s+1} -h), \ldots, \varepsilon_k \alpha_k) &  \mbox{ if } \varepsilon_{s+1} = \omega^j, \\
(\varepsilon_1 \alpha_1, \ldots, \varepsilon_{s+1} h, \omega^j 1, \varepsilon_{s+1}(\alpha_{s+1}-h), \ldots, \varepsilon_k \alpha_k) & \mbox{ otherwise.}
\end{cases}
\]
These cases are again distinct and provide the remaining covers for $\C^{(m)}(u)$.
\end{proof}

Theorem \ref{thm:saturated} yields several easy corollaries. The first is analogous to property Y2 of Young's lattice; the second corresponds to Y3.

\begin{cor}\label{cor:f}
The number of saturated chains from $\emptyset$ to $\alpha$ in $\co^{(m)}$ is equal to the number $f_n^{(m)}(\alpha)$ of $m$-colored permutations $w$ with colored descent composition $\alpha$.
\end{cor}

\begin{cor}
The total number of saturated chains from $\emptyset$ to rank $n$ is equal to the number of $m$-colored permutations of $[n]$, \[\sum_{\alpha \in \com^{(m)}(n)} f_n^{(m)}(\alpha) =  m^n n!.\]
\end{cor}

\begin{cor}
The number of $m$-colored compositions $\beta \in \com^{(m)}(n+1)$ covering $\alpha \in \com^{(m)}(n)$ is $m(n+1)$.
\end{cor}


\section{Colored quasisymmetric functions}\label{sec:quasi}

One key use for compositions is as an indexing set for quasisymmetric functions. Similarly, there exist colored quasisymmetric functions (due to Poirier \cite{Poirier}) that use colored compositions as indices. Both these situations are analogous to how partitions index symmetric functions.

Recall (\cite{Stanley2}, ch. 7.19) that a quasisymmetric function is a formal series \[Q(x_1, x_2, \ldots ) \in \mathbb{Z}[[x_1, x_2,\ldots ]] \] of bounded degree such that for any composition $\alpha = (\alpha_1, \alpha_2, \ldots, \alpha_k)$, the coefficient of $x_{i_1}^{\alpha_1} x_{i_2}^{\alpha_2} \cdots x_{i_k}^{\alpha_k}$ with $i_1 < i_2 < \cdots < i_k$ is the same as the coefficient of $x_1^{\alpha_1} x_2^{\alpha_2} \cdots x_k^{\alpha_k}$. One natural basis for the quasisymmetric functions homogeneous of degree $n$ is given by the \emph{fundamental quasisymmetric functions}, $L_{\alpha}$, where $\alpha$ ranges over all of $\com(n)$. If $\alpha = (\alpha_1, \ldots, \alpha_k) \models n$, then define
\[ L_{\alpha} := \sum x_{i_1} \cdots x_{i_n},\] where the sum is taken over all $i_1 \leq i_2 \leq \cdots \leq i_n$ with $i_s < i_{s+1}$ if $s = \alpha_1 + \cdots + \alpha_r$ for some $r$. For example, \[ L_{21} = \sum_{i \leq j < k} x_i x_j x_k.\]

Colored quasisymmetric functions are simply a generalization of quasisymmetric functions to an alphabet with several colors for its letters. For fixed $m$, we consider formal series in the alphabet \[ X^{(m)} := \{ x_{0,1}, x_{0,2}, \ldots, x_{1,1}, x_{1,2}, \ldots, x_{m-1,1}, x_{m-1,2}, \ldots \},\] (so the first subscript corresponds to color) with the same quasisymmetric property. Namely, an $m$-colored quasisymmetric function $Q(X^{(m)})$ is a formal series of bounded degree such that for any $m$-colored composition $\alpha= (\omega^{j_1} \alpha_1, \ldots, \omega^{j_k} \alpha_k)$, the coefficient of $x_{j_1,i_1}^{\alpha_1} x_{j_2,i_2}^{\alpha_2} \cdots x_{j_k,i_k}^{\alpha_k}$ with $i_1 < i_2 < \cdots < i_k$ is the same as the coefficient of $x_{j_1,1}^{\alpha_1} x_{j_2,2}^{\alpha_2} \cdots x_{j_k,k}^{\alpha_k}$. Intuitively, the letters are colored the same as the parts of $\alpha$. The \emph{$m$-colored fundamental quasisymmetric functions} are defined as follows. First, if $s = \alpha_1 + \cdots + \alpha_r + h$, $1 \leq h \leq \alpha_{r+1}$, then define $j'_s = j_{r+1}$, the color of part $\alpha_{r+1}$. Then,
\[L^{(m)}_{\alpha} := \sum x_{j'_1,i_1} \cdots x_{j'_n,i_n}, \] where the sum is taken over all $i_1 \leq i_2 \leq \cdots \leq i_n$ with $i_s < i_{s+1}$ if both $j'_s \geq j'_{s+1}$ and $s = \alpha_1 + \cdots + \alpha_r$ for some $r$.  For example, \[ L^{(2)}_{1\bar{2}\bar{1}} = \sum_{i \leq j \leq k < l} x_i y_j y_k y_l \quad \mbox{ and } \quad L^{(3)}_{2\bar{\bar{ 1}} \bar{2}} = \sum_{i \leq j \leq k < l \leq m} x_i x_j z_k y_l y_m.\] As in the ordinary case, the $L^{(m)}_{\alpha}$, where $\alpha$ ranges over $\com^{(m)}(n)$, give a basis for the $m$-colored quasisymmetric functions homogeneous of degree $n$.

There is a nice formula for multiplying colored quasisymmetric functions in the fundamental basis. Let $u \in C_m \wr \mathfrak{S}_n$ and let $v$ be an $m$-colored permutation of the set $\{ n+1, n+2, \ldots, n+r\}$. Let $\alpha = \C^{(m)}(u)$ and $\beta = \C^{(m)}(v)$. Then we have \[ L^{(m)}_{\alpha} L^{(m)}_{\beta} = \sum_{w} L^{(m)}_{\C^{(m)}(w)},\] where the sum is taken over all \emph{shuffles} $w$ of $u$ and $v$, i.e., all colored permutations $w \in C_m \wr \mathfrak{S}_{n+r}$ such that $w[n] = u$ and $w$ restricted to $\{ n+1, n+2, \ldots, n+r\}$ is $v$.

If $r =1$, then we see that the shuffles of $u$ and $v = \omega^j(n+1)$ are precisely those permutations $u_{(i,j)}$ from the proof of Theorem \ref{thm:saturated}. Applying the multiplication rule, and summing over all $j$, we have a Pieri-type rule analogous to property Y4 of Young's lattice.

\begin{prp}\label{prp:pieri}
We have: \[(L^{(m)}_1 + L^{(m)}_{\omega 1} + \cdots + L^{(m)}_{\omega^{m-1} 1}) L^{(m)}_{\alpha} = \sum_{ \alpha \prec \beta } L^{(m)}_{\beta}. \]
\end{prp}

As Bj\"{o}rner and Stanley remark in the case of a single color, we
could have used Proposition \ref{prp:pieri} to define the poset
$\co^{(m)}$ in the first place. At the least, it is a good
justification for the study of $\co^{(m)}$.

Repeated application of the proposition gives the formula
\[(L^{(m)}_1 + L^{(m)}_{\omega 1} + \cdots + L^{(m)}_{\omega^{m-1} 1})^n = \sum_{\alpha \in \com^{(m)}(n)} f_n^{(m)}(\alpha) L^{(m)}_{\alpha},\] where $f_n^{(m)}$ is the number of $m$-colored permutations with colored descent composition $\alpha$. This equation is equivalent to Corollary \ref{cor:f}, and analogous to the following formula for Schur functions (see \cite{Stanley2}) that corresponds to property Y2 of Young's lattice:
\[ s_1^n = \sum_{\lambda \vdash n} f^{\lambda} s_{\lambda},\] where $f^{\lambda}$ is the number of Young tableaux of shape $\lambda$.

\section{Shellability and M\"{o}bius function}\label{sec:mobius}
In this section we show that $\co^{(m)}$ is CL-shellable by giving
an explicit dual CL-labeling.  See \cite{BGS} for an introduction to
CL-shellable posets.  We use a model of removing colored balls from
urns to define our labeling on downward maximal chains. Given a
colored composition of $n$, $\alpha = (\varepsilon_1 \alpha_1,
\ldots, \varepsilon_k \alpha_k)$, we picture $k$ urns next to each
other, labeled $U_1, U_2, \ldots, U_k$ from left to right. In urn
$U_i$ we start with $\alpha_i$ balls of color $\varepsilon_i$, for a
total of $n$ balls. Moving down along a maximal chain,
we remove a ball from an urn for each covering relation, and possibly move
some balls from one urn to another. There are three different types
of moves, which we now describe. After some number of steps, suppose
that $U_i$ is a nonempty urn, $U_h$ is the first nonempty urn on its
left, and $U_j$ is the first nonempty urn on its right. Let
$\beta_i, \beta_h, \beta_j$ be the number of balls in the
corresponding urns and let $\varepsilon_i, \varepsilon_h,
\varepsilon_j$ be the colors of those balls.  The three possible
moves are:
\begin{enumerate}
\item If $\beta_i \geq 2$, or if $\varepsilon_h \neq \varepsilon_i$, or if $U_i$ is the first nonempty urn, then remove a ball from urn $U_i.$

\item If $\beta_i = 1$ and $\varepsilon_h = \varepsilon_j \neq \varepsilon_i$, then remove the ball from $U_i$ and place all the balls from $U_h$ and $U_j$
into $U_i.$

\item If $\beta_i \geq 2$ and $\varepsilon_j = \varepsilon_i$, then move all balls from $U_j$ to $U_i$ and remove a ball from
$U_i.$

\end{enumerate}
After any number of moves, we may associate the distribution of
colored balls in urns with an element of $\co^{(m)}$. The
different urns represent the parts of the composition, the number of
balls in an urn is the size of that part, and the color of the balls
is the color of the part.  Here we ignore parts of size 0.  Notice
that the color of a part is well defined, since none of the moves
allows balls of different colors to be combined in a single urn.  It
is an easy exercise to check that each covering relation in
$\co^{(m)}$ corresponds to one of these three possible moves for
some urn, and furthermore that the urn and type of move are unique.

Let $[\emptyset, \alpha]$ be an interval in $\co^{(m)}$, $|\alpha| =
n.$  We will now define a labeling \[\lambda(c)= (\lambda_1(c),
\lambda_2(c), \ldots, \lambda_n(c))\] for a maximal chain \[ c = (
\alpha = \alpha^0 \succ \alpha^1 \succ \cdots \succ \alpha^n =
\emptyset ).\] Our set of labels is $\mathbb{N}\times \{1,2,3\},$
totally ordered with the lexicographic order.  For each covering
relation $\alpha^{r-1} \succ \alpha^{r}$ we have a unique urn and
type of move that takes the distribution of balls in urns for
$\alpha^{r-1}$ to the distribution for $\alpha^{r}.$ Suppose that
move is of type $t,$ and removes a ball from urn $U_i.$ Then we
define the label $\lambda_{r}(c) = (i, t).$

Notice that with labels defined on maximal chains in lower intervals
$[\emptyset, \alpha],$ there is an induced labeling defined on
maximal chains in arbitrary intervals $[\beta, \alpha].$  As an
example, consider the following two maximal chains in $[3,
22\bar{1}2]$:
\begin{align*}
c_0 & = ( 22\bar{1}2 \succ 12\bar{1}2 \succ 2\bar{1}2 \succ 1\bar{1}2 \succ 3) \\
c & = ( 22\bar{1}2 \succ 21\bar{1}2 \succ 2\bar{1}2 \succ 22 \succ
3)
\end{align*}
They are labeled $\lambda(c_0) = ((1,1),(1,1),(2,1),(3,2))$ and
$\lambda(c) = ((2,1),(1,3),(3,1),(1,3))$. Pictured as colored balls
and urns, we have:
\begin{align*}
c_0 & :  \lfloor \bullet \bullet\rfloor \lfloor\bullet
\bullet\rfloor \lfloor\circ\rfloor \lfloor\bullet \bullet\rfloor
\overset{}{\rightarrow} \lfloor \bullet \rfloor \lfloor\bullet
\bullet\rfloor \lfloor\circ\rfloor \lfloor\bullet \bullet\rfloor
\overset{}{\rightarrow} \lfloor \rfloor \lfloor\bullet
\bullet\rfloor \lfloor\circ\rfloor \lfloor\bullet \bullet\rfloor
\overset{}{\rightarrow} \lfloor\rfloor \lfloor \bullet\rfloor
\lfloor\circ\rfloor \lfloor\bullet \bullet\rfloor
\overset{}{\rightarrow} \lfloor \rfloor \lfloor\rfloor
\lfloor\bullet\bullet\bullet\rfloor \lfloor
\rfloor \\
c &:  \lfloor \bullet \bullet\rfloor \lfloor\bullet \bullet\rfloor
\lfloor\circ\rfloor \lfloor\bullet \bullet\rfloor
\overset{}{\rightarrow} \lfloor \bullet \bullet \rfloor
\lfloor\bullet\rfloor \lfloor\circ\rfloor \lfloor\bullet
\bullet\rfloor \overset{}{\rightarrow} \lfloor \bullet\bullet\rfloor
\lfloor\rfloor \lfloor\circ\rfloor \lfloor\bullet \bullet\rfloor
\overset{}{\rightarrow} \lfloor\bullet\bullet\rfloor \lfloor \rfloor
\lfloor\rfloor \lfloor\bullet \bullet\rfloor \overset{}{\rightarrow}
\lfloor \bullet\bullet\bullet\rfloor \lfloor\rfloor \lfloor\rfloor
\lfloor \rfloor .
\end{align*}
In fact, the chain $c_0$ above is lexicographically minimal and has the only increasing label. Notice that if we start with all balls of the same color, moves of type $(2)$
cannot occur and we recover the construction in the appendix of
\cite{BjornerStanley}.  These labels agree with the labels defined there, with $(i, 1) \mapsto i$
and $(i, 3) \mapsto i^{\prime}$.

By proving that this labeling is in fact a CL-labeling, we obtain
our analog of property Y5 of Young's lattice.  The proof is given
in section \ref{sec:clproof}.

\begin{thm}\label{thm:shellable}
Intervals in $\co^{(m)}$ are dual CL-shellable and hence
Cohen-Macaulay.
\end{thm}


Now we calculate the M\"{o}bius function of lower intervals. As
always, we must have $\mu_{\co^{(m)}} (\emptyset, \emptyset)=1.$ For
$\alpha \neq \emptyset,$ we have the following.

\begin{prp}
$$ \mu_{\co^{(m)}} (\emptyset, \alpha) = \begin{cases}
(-1)^{|\alpha|} & \mbox{if } \alpha = (\varepsilon_1 1,
\varepsilon_2 1, \ldots, \varepsilon_{|\alpha| } 1) \\ & \mbox{for
some colors } \varepsilon_1 \neq \varepsilon_2 \neq \cdots \neq
\varepsilon_{|\alpha|}, \\  0 & \mbox{otherwise.}
\end{cases}$$
\end{prp}

\begin{proof}
We make use of the combinatorial description of the M\"{o}bius
function for a graded poset with a CL-labeling, given in \cite{BGS}.
That is, the M\"{o}bius function of an interval is $-1$ to the
length of the interval, times the number of maximal chains with a
strictly decreasing label.

Suppose that $\alpha = (\varepsilon_1 1, \varepsilon_2 1, \ldots,
\varepsilon_{|\alpha| } 1),$ with $\varepsilon_1 \neq \varepsilon_2
\neq \cdots \neq \varepsilon_{|\alpha|}.$  Then there is a unique
chain with a strictly decreasing label, obtained by removing the
balls from right to left using only type (1) moves.  Therefore
$\mu_{\co^{(m)}} (\emptyset, \alpha) = (-1)^{|\alpha|}.$

Now suppose that $\alpha$ has a part $i$ of size 2 or greater.  We
want to show that there is no chain in $[\emptyset, \alpha]$ with a
strictly decreasing label.  Any chain that makes a type (1) move
from the same urn twice will have a repeated label.  The only way to
remove the balls from urn $i$ and possibly have a decreasing label
is to remove at most one ball with a type (1) move, and then move
all the balls to an urn on the left with a type (2) or (3) move. But
in the new urn we have at least two balls, and the process repeats.
At some point we must have an urn with at least two balls and no way
to make a type (2) or (3) move.  Then we must use two type (1) moves
for the same urn, so the chain label cannot be strictly decreasing.

Finally, suppose that $\alpha$ has parts $\alpha_i$ and
$\alpha_{i+1}$ of size 1 and the same color.  The only legal way to
remove the balls from the corresponding urns is to remove the left
one first, creating an increase in the chain label.

\end{proof}

Note that for $\alpha \models_m n$ with $\mu_{\co^{(m)}} (\emptyset,
\alpha) \neq 0,$ there are $m$ choices for the color of the first
part, and $m-1$ choices for the color of each succeeding part. Hence
there are $m(m-1)^{n-1}$ compositions $\alpha \models_m n$ with
$\mu_{\co^{(m)}} (\emptyset, \alpha) \neq 0.$  For $m > 1$ an
elementary calculation gives the following generating function. \[
\displaystyle\sum_{\alpha \in \co^{(m)}} \mu_{\co^{(m)}} (\emptyset,
\alpha) t^{|\alpha |} = \frac{1+t}{1-(m-1)t}.\]

Define the following ``truncated" poset, \[\co^{(m)}_n := \widehat{1} \cup \bigcup_{1 \leq i
\leq n} \com^{(m)}(i),\] with the order relation as before except
with a new maximal element $\widehat{1}$ that covers all the
compositions in $\com^{(m)}(n)$.

\begin{cor}
The poset $\co^{(m)}_n$ is shellable, with M\"{o}bius
function
$$\mu(\emptyset, \widehat{1}) = (-1)^{n+1} (m-1)^n.$$
\end{cor}

The proof of this corollary follows the argument of
\cite{BjornerStanley}.

\begin{proof}

First we want to show that every $m$-colored composition $\alpha =
(\varepsilon_1 \alpha_1, \varepsilon_2 \alpha_2, \ldots,
\varepsilon_k \alpha_k)$ of at most $n$ lies below the composition
$\gamma^n \models_m mn$, defined as the concatenation of $n$ copies
of $\gamma = (1, \omega 1, \omega^2 1, \ldots, \omega^{m-1} 1)$. To
the $i^{\rm th}$
 part of $\alpha$ we can associate  $\alpha_i$
copies of $\gamma$.  First, we use covering relations of type $(2)$
(as originally described), $\alpha_i$ times to split the part into
all parts of size 1 and color $\varepsilon_i$.  Then we use covering
relations of type $(3)$ to fill in the remaining 1's of different
colors. Therefore $\co^{(m)}_n$ is obtained via rank
selection from the interval $[\emptyset, \gamma^n],$ so shellability
follows by results of \cite{BjornerShellable}.

For the M\"{o}bius function:

\begin{align*}
\mu(\emptyset, \widehat{1}) = - \sum_{ |\alpha| \leq n}
\mu(\emptyset, \alpha) & = -\left(1 + \sum_{k=1}^n (-1)^k
m(m-1)^{k-1}\right)
\\ & = -\left(1 - m \sum_{k=0}^{n-1} (1-m)^{k}\right) \\ & = -\left( 1 -
m \left( \frac{1-(1-m)^n}{1-(1-m)}\right) \right) \\ &=
(-1)^{n+1}(m-1)^n.
\end{align*}

\end{proof}

\begin{remark}\label{rem:1}  Bj\"{o}rner and Stanley show that
$\co^{(1)}$ is isomorphic to the subword order on 2 letters. This
allows the transfer of many results, such as CL-shellability and the
M\"{o}bius function of an arbitrary interval.  Since there are
$m(m+1)^{n-1}$ colored compositions of $n,$ one might expect that
$\co^{(m)}$ is isomorphic to the subword order on $m+1$ letters,
with a restriction on the first or last letter of a word.  However,
this turns out to be false.  For example, consider the colored
composition $(1, \bar{1}, 1)$, which is present in $\co^{(m)}$ for all $m \geq 2$.  It covers 4 colored compositions: $(2), (1,\bar{1}), (\bar{1},1), (1,1)$,
but a word of length 3 can cover at most 3 subwords.  Therefore the
results of this section do not follow directly from properties of
the subword order.
\end{remark}


\section{Proof of CL-shellability}\label{sec:clproof}
In this section we give the proof of Theorem \ref{thm:shellable}. We first note that the labeling is a well defined chain labeling. That is, if two chains agree on their first $k$ edges, then their
first $k$ labels agree.  This is clear from the definition.

The labeling of maximal chains in $[\emptyset, \alpha]$ gives an
induced labeling on rooted intervals $([\beta, \alpha], c),$ where
$c$ is a maximal chain in $[\emptyset, \beta].$  This induced
labeling is of the same kind, so it suffices to check that the
properties of a CL-labeling hold in an arbitrary interval $[\beta,
\alpha].$

First, we want to show that the chain with the lexicographically
first label has a weakly increasing label. The lexicographically
first label is well defined, since all the moves from a given
distribution of balls in urns have distinct labels. Moreover, we can describe the lexicographically first chain, $c_0$, as follows. If \[ c_0 = ( \alpha = \alpha^0 \succ \alpha^1 \succ \cdots \succ \alpha^k = \beta),\] then at each step, to move down from $\alpha^{r-1}$ to $\alpha^r$, we must remove a ball from an urn as far to the left as
possible, such that the new composition is still in the interval $[\beta, \alpha]$. To prove that $\lambda_1(c_0) \leq \lambda_2(c_0) \leq \cdots \leq \lambda_k(c_0)$, we use the following lemma.

\begin{lem}\label{lem:1} On an interval of length two, the chain with the
lexicographically first label is weakly increasing.
\end{lem}

\begin{proof}
On an interval of length two, all chains correspond to removing two
balls from urns, such that the starting and ending distributions are
the same.  For the chain $c_0$ with the lexicographically first
label, the urns are consecutively chosen to be as far to the left as
possible. Suppose that one ball is removed from urn $U_i$
and the other ball is removed from urn $U_j.$

If $j > i+1,$ then there is a nonempty urn between $U_i$ and $U_j,$
and removing a ball from one of the urns does not affect the
possibility of removing the other ball from its urn.  Therefore it
is clear that $\lambda(c_0)$ is weakly increasing.

Suppose that $j = i+1,$ so there is no urn between $U_i$ and $U_j.$
Removing a ball from $U_i$ does not affect the possibility of
removing a ball from $U_{i+1},$ unless $\varepsilon_i \neq
\varepsilon_{i+1} = \varepsilon_{i-1}$ and $\alpha_i = 1.$  In this
case, the urns could have been chosen to be $U_h$ and $U_i,$ for an
appropriate urn $U_h$ with $h < i,$ making a type (1) or type (3)
move in $U_h$ and then a type (2) move in urn $U_i.$ However, this
contradicts our assumption that $U_i$ was chosen to be the leftmost
possible. So if $i \neq j,$ $\lambda(c_0)$ is weakly increasing.

The only remaining case is if $U_i = U_j.$  If $\alpha_i > 2,$ then
we have the weakly increasing label $\lambda(c_0) = ((i, 1), (i,
1)).$  Now suppose that $\alpha_i = 2$ and $\varepsilon_{i-1} \neq
\varepsilon_i.$  Then we also have $\lambda(c_0) = ((i, 1), (i,
1)).$  If $\varepsilon_{i-1} = \varepsilon_i,$ then $c_0$ is found
by choosing an appropriate urn $U_h, h < i$ and making a type (1) or
(3) move in urn $U_h$ and then a type (1) move in urn $U_i.$  Again,
this contradicts our assumption that $U_i$ was chosen to be the
leftmost possible. Therefore if $i = j,$ $\lambda(c_0)$ is weakly
increasing.
\end{proof}

Returning to the general case, for every $r$, the induced labeling of $c_0$ on the chain $\alpha^{r-1}
\succ \alpha^r \succ \alpha^{r+1}$ is lexicographically first on the
interval $[\alpha^{r+1}, \alpha^{r-1}]$. Then by Lemma \ref{lem:1}, $\lambda_1(c_0) \leq \lambda_2(c_0) \leq \cdots \leq \lambda_k(c_0)$, i.e., $\lambda(c_0)$ is weakly increasing.
Now it remains only to show no chain other than $c_0$ has a weakly
increasing label.

If another chain results with the same distribution of balls into
urns as the lexicographically first chain, including the locations of the empty
urns, then the label on that chain must have a descent. To see this,
consider the point where it deviates from the lexicographically first chain. It is
leaving a ball behind in a lower numbered urn.  At some later step
it must remove a ball from that urn, which will create a descent.

Now we need to consider chains which result in the same distribution
of balls into nonempty urns as the lexicographically first chain, but such that the empty
urns are in different positions.  Let $c$ be such a chain.  Since
the lexicographically first chain removes balls from urns from left to right, the
final distribution of balls into urns for the lexicographically first chain has
its nonempty urns as far to the right as possible.

Let $U_1, U_2, \ldots$ be the urns labeled from left to right.  Let
$U_i(c)$ be the number of balls in urn $U_i$ in the final
distribution for the chain $c$.  Let $r$ be the largest number such
that $U_r(c) > 0$ and $U_r(c_0) = 0.$  There is such an $r$ by our
assumption on $c.$  We must have a $j$ such that: 1) $U_{r+i}(c) = 0$
for all $1 \leq i \leq j$, 2) $U_{r+i}(c_0) = 0$ for all $1 \leq i
< j$, and 3) $U_{r+j}(c_0) \neq 0.$  Note that the color and number of
the balls in urn $U_r$ for $c$ is the same as the color and number
of the balls in urn $U_{r+j}$ for $c_0.$

If $c$ has an increasing label, then the urns $U_{r+1}, \ldots,
U_{r+j}$ must be emptied from left to right.  Therefore at some
point the urns $U_{r+1}, \ldots, U_{r+j-1}$ are empty and we need to
remove the last ball from $U_{r+j}.$  The only way to do this is to
use a move of type (2), which has a label $(r, 2)$ and creates a
descent.


\begin{thebibliography}{99}

\bibitem{BjornerShellable} A. Bj\"{o}rner, \emph{Shellable and
Cohen-Macaulay posets}, Trans. Amer. Math. Soc. {\bf 260} (1980),
159--183.


\bibitem{BGS} A. Bj\"{o}rner, A. Garsia and R. Stanley,
\emph{Cohen-Macaulay partially ordered sets} in \emph{Ordered Sets},
(I. Rival, ed.), Reidel, Dordrecht/Boston, 1982, pp. 583--615.

\bibitem{BjornerStanley} A. Bj\"{o}rner and R. Stanley, \emph{An analogue of Young's lattice for compositions}, arXiv: math.CO/0508043.

%

\bibitem{Poirier} S. Poirier, \emph{Cycle type and descent set in wreath products}, Discrete Mathematics {\bf 180} (1998), 315--343.


\bibitem{Stanley2} R. Stanley, \emph{Enumerative Combinatorics, Volume II}, Cambridge University Press, 2001.


\end{thebibliography}
\end{document}